\documentclass[10pt]{article}
\usepackage{amsmath,amssymb}

\newcommand{\hess}{\operatorname{Hess}}
\newcommand{\grad}{\operatorname{grad}}
\newcommand{\isom}{\operatorname{Isom}}

\begin{document}

\subsection*{WP Metric Geometry Quick View \quad\quad SAW  July 3, 2007}

Let $\mathcal T$ be the Teichm\"{u}ller space for marked genus $g$, $n$ punctured Riemann surfaces $R$ of negative Euler characteristic.  By Uniformization a conformal structure determines a unique complete compatible hyperbolic metric $ds^2$.  $\mathcal T$ is a complex manifold of $\dim_{\mathbb C}\mathcal=3g-3+n$ with the cotangent space at $R$ represented by $Q(R)$, the space of holomorphic quadratic differentials on $R$ with at most simple poles at the punctures.  The Weil-Petersson Hermitian cometric is 
\[
\langle\varphi,\psi\rangle=\int_R\varphi\overline{\psi}(ds^2)^{-1}.
\]
The WP metric is invariant under the action of the {\em mapping class group} MCG, \cite{Ivmcg}.  The metric is K\"{a}hler, non complete with negative sectional curvature with $\sup_{\mathcal T}=0$ (except for $\dim \mathcal T=1$) and $\inf_{\mathcal T}=-\infty$.  
{\em In practice and experience the WP geometry of $\mathcal T$ corresponds to the hyperbolic geometry of surfaces.}

{\bf The augmented Teichm\"{u}ller space,} \cite{Abdegn, Bersdeg}.  The augmented Teichm\"{u}ller space $\overline{\mathcal T}$ is a MCG invariant bordification (a partial compactification) following the approach for Baily-Borel bordifications. ($\overline{\mathcal T}$ is an analog of the $SL(2;\mathbb Z)$ invariant bordification $\mathbb Q$  for $\mathbb H$ provided with the Satake horocycle topology.)  The bordification is also described as the Chabauty topology closure of the faithful cofinite representations of $\pi_1(R)$ into $PSL(2;\mathbb R)$.  $\overline{\mathcal T}/MCG$ is homeomorphic to the Deligne-Mumford stable curve compactification of the moduli space of Riemann surfaces.
The elements of $\overline{\mathcal T}-\mathcal T$ are marked degenerate hyperbolic structures for which certain simple non peripheral closed curves are {\em represented} by pairs of cusps. 

{\bf CAT(0) geometry.}  Important is that $\overline{\mathcal T}$ is WP complete with 
$(\overline{\mathcal T} ,d_{WP})$ a $CAT(0)$ metric space (a generalized complete, simply connected, non positively curved space.)   The large-scale geometry is described by $\overline{\mathcal T}$ being quasi isometric to $\mathcal {PG}(R)$ the {\em pants graph} of $R$ with unit-length edge metric, \cite{Brkwp}.  $\overline{\mathcal T}$ is a stratified space with each open strata characterized as the union of all geodesics (distance realizing paths) containing a given point as an interior point. $\overline{\mathcal T}$ itself is characterized as the closed WP convex hull of the maximally degenerate hyperbolic structures (the unions of thrice punctured spheres.)  An application of the $CAT(0)$ geometry and the important rigidity of the {\em complex of curves} $\mathcal C(R)$, \cite[Chap. 3]{Ivmcg} is that the isometry group $\isom_{WP}$ coincides with MCG. 

{\bf Geodesic-length functions.}  Associated to each non trivial, non peripheral free homotopy class on $R$ is the length of the unique $ds^2$ geodesic in the homotopy class.  Geodesic-lengths $\ell_{\alpha}(R)$ on $\mathcal T$ are elements of the WP geometry.  The Fenchel-Nielsen twist (right earthquake) deformation about $\alpha$ is described by $2t_{\alpha}=J\grad\ell_{\alpha}$ for $J$ the Teichm\"{u}ller almost complex structure.  The WP Hermitian pairing is described in terms of the hyperbolic trigonometry for $R=\mathbb H/\Gamma$ with
\[
\langle\grad\ell_{\alpha},J\grad\ell_{\beta}\rangle=4\,\omega_{WP}(t_{\alpha},t_{\beta})=-2\sum_{p\in\alpha\cap\beta}\cos\theta_p
\]
for the geodesics $\alpha,\,\beta$ and intersection angles $\theta_*$ on $R$ and in upper half plane  
\[
\langle\grad\ell_{\alpha},\grad\ell_{\beta}\rangle=\frac{2}{\pi}\bigl(\ell_{\alpha}\delta_{\alpha\beta}\ +{\ 
\sum}^{\prime}_{\langle A\rangle\backslash\Gamma/\langle B\rangle}(u\log\frac{u+1}{u-1}-2\bigr)\bigr)
\]
for the Kronecker delta $\delta_*$, where for $C\in\langle A\rangle\backslash\Gamma/\langle B\rangle$ then $u=u(\tilde\alpha,C(\tilde\beta))$ is the cosine of the intersection angle if the lifts $\tilde\alpha$ and $C(\tilde\beta)$ intersect and is otherwise $\cosh d(\tilde\alpha,C(\tilde\beta))$.

{\bf WP convexity and curvature.}  Information on geodesics is provided by the strict convexity of geodesic-lengths (more generally convexity of the total length of a {\em measured geodesic lamination}) along WP geodesics (square roots are also convex.) For surfaces with cusps the distances along geodesics between horocycles are also strictly convex.  Geodesic-length sublevel sets are convex.  On $\mathcal T$ the WP Levi-Civita connection $D$ is described for small geodesic-length, $\lambda=\grad\ell^{1/2}_{\alpha}$ and tangent $U$ by
\[
D_U\lambda_{\alpha}\ =\ 3\ell_{\alpha}^{-1/2}\langle J\lambda_{\alpha},U\rangle J\lambda_{\alpha}\ +\ O(\ell_{\alpha}^{3/2}\|U\|).
\]
The remainder term constant is uniform for bounded geodesic-length.  Bounds for the gradient and Hessian of geodesic-length give rise to applications.  The WP curvature of the span $\{t_{\alpha},Jt_{\alpha}\}$ is $O(-\ell_{\alpha}^{-1})$.  Similarly for a pair of deformations supported on different components of $R-\{short\ geodesics\}$ the curvature of the corresponding $2$-plane is $O(-short\ separating\ lengths)$.

{\bf Fenchel-Nielsen coordinates.}  Marked hyperbolic pairs of pants are determined by their three boundary geodesic-lengths in $\mathbb R_+$.  The FN twist-length coordinates $(\ell_j,\tau_j)^{3g-3+n}_{j=1}$ for assembling hyperbolic pants provide global coordinates for $\mathcal T$ with expansions
\[
\omega_{WP}=\frac12\sum\,d\ell_j\wedge d\tau_j
\]
and on the Bers region $\{\ell_j<c_0\}$ the expansions 
\[
\langle\ ,\ \rangle \ \asymp\ \sum  (d\ell_j^{1/2})^2+(d\ell_j^{1/2}\circ J)^2\ \asymp\ 
\sum \hess\ell_j 
\]
with comparability uniform in terms of $c_0$.   At the maximally degenerate structure the WP metric has the expansions 
\begin{align*}
\langle\ ,\ \rangle \ =&\ 2\pi\sum (d\ell_j^{1/2})^2+(d\ell_j^{1/2}\circ J)^2\ +\ O(\sum \ell_j^3\,\langle\ ,\ \rangle )\\ =&\ \frac{\pi}{6}\sum \frac{\hess\ell_j^2}{\ell_j} \ +\ O(\sum \ell_j^2\,\langle\ ,\ \rangle ).
\end{align*}

{\bf Alexandrov tangent cone.}  For a pair of geodesics (distance realizing curves) from a common point in a $CAT(0)$ metric space there is a well-defined initial angle, \cite[Chap. II.3]{BH}.  A tangent cone is defined in terms of initial angle.  At a point $p$ of a strata $\mathcal T(\sigma)\subset \overline{\mathcal T}-\mathcal T$, $\sigma\in\mathcal C(R)$, the WP Alexandrov tangent cone $AC_p$ is isometric to a product of a Euclidean orthant and the tangent space $\mathbf T\mathcal T(\sigma)$ with WP pairing.  The isometry is given in terms of initial derivatives of geodesic-length functions.  The dimension of the Euclidean orthant is the count of geodesic-lengths trivial on $\mathcal T(\sigma)$ (the count of nodes.)   A property of $\overline{\mathcal T}$ is that a WP geodesic tangential to $\mathcal T(\sigma)$ at $p$ actually lies in $\mathcal T(\sigma)$.  Applications include that geodesics do not {\em refract} at the bordification and that the angle of incidence and reflection coincide for certain limits of degenerating geodesics.  A further application is for combinatorial harmonic maps.  Certain groups acting on Euclidean buildings and group extensions acting on Cayley graphs satisfying a Poincar\'{e} type inequality for links of points will have a global fixed point for an action on $\overline{\mathcal T}$.

All important attributions and references are provided in the introductions of \cite{Wlcomp,Wlbhv}.   
  

\end{document}